\documentclass[a4paper]{amsart}
\pdfoutput=1

\usepackage{amssymb, amsmath}
\usepackage[english]{babel}
\usepackage{ifthen}
\usepackage{ifpdf}
\usepackage{graphicx}
\usepackage{tikz}
\usepackage{url}
\usepackage{caption}
\usepackage{amscd}
\usepackage{cite}

\usepackage[pdftex, plainpages = false, pdfstartview=FitH,
    pdfpagelabels,  pdfpagelayout = SinglePage, 
    bookmarks = true, bookmarksopen = true, 
    bookmarksnumbered = true, linktocpage,
    colorlinks = true, 
    linkcolor = blue, urlcolor  = blue, 
    citecolor = red,  anchorcolor = green, 
    hyperindex = true, hyperfigures]{hyperref} 

\DeclareGraphicsExtensions{.png, .jpg, .pdf}
\graphicspath{{img/}{img/}{img/}}

\title[Toric line arrangements]%
			{Face enumeration for line arrangements in a $2$-torus}
\author[K\, Chandrashekhar]{Karthik Chandrashekhar}
\author[P\, Deshpande ]{Priyavrat Deshpande}
\address{Chennai Mathematical Institute \\ Chennai\\ India}
\email{aneesh@cmi.ac.in}
\email{pdeshpande@cmi.ac.in} 


\numberwithin{equation}{section}

\theoremstyle{plain}
\newtheorem{theorem}{Theorem}[section]
\newtheorem{lemma}[theorem]{Lemma}
\newtheorem{cor}[theorem]{Corollary}
\newtheorem{proposition}[theorem]{Proposition}

\theoremstyle{definition}
\newtheorem{defn}[theorem]{Definition}
\newtheorem{example}[theorem]{Example}

\theoremstyle{remark}

\newcommand{\bt}[1]{\begin{theorem}\label{#1}}
\newcommand{\bc}[1]{\begin{cor}\label{#1}}
\newcommand{\bl}[1]{\begin{lemma}\label{#1}}
\newcommand{\bp}[1]{\begin{proposition}\label{#1}}
\newcommand{\be}[1]{\begin{example}\label{#1}}
\newcommand{\bd}[1]{\begin{defn}\label{#1}}
\newcommand{\br}[1]{\begin{remark}\label{#1}}
\newcommand{\bx}[1]{\begin{exercise}\label{#1}}
\newcommand{\bcon}[1]{\begin{conjecture}\label{#1}}

\newcommand{\et}{\end{theorem}}
\newcommand{\ec}{\end{cor}}
\newcommand{\el}{\end{lemma}}
\newcommand{\ep}{\end{proposition}}
\newcommand{\ee}{\end{example}}
\newcommand{\ed}{\end{defn}}
\newcommand{\exc}{\end{exercise}}
\newcommand{\er}{\end{remark}}
\newcommand{\econ}{\end{conjecture}}

\newcommand{\bpr}{\begin{proof}}
\newcommand{\epr}{\end{proof}}

\newcommand{\ol}{\overline}

\def\A  {\mathcal{A}}

\def \F {\mathcal{F}}

\def \h {\mathbb{H}}

\def\R{\mathbb{R}}
\def\N{\mathbb{N}}

\def\Z{\mathbb{Z}}

\def\t{\mathbb{T}}

\def\ds{\displaystyle}

\begin{document}

\begin{abstract}
A toric arrangement is a finite collection of codimension-$1$ subtori in a torus. These subtori stratify the ambient torus into faces of various dimensions. Let $f_i$ denote the number of $i$-dimensional faces; these so-called face numbers satisfy the Euler relation $\sum_i (-1)^i f_i = 0$. However not all tuples of natural numbers satisfying this relation arise as face numbers of some toric arrangement. In this paper we focus on toric arrangements in a $2$-dimensional torus and obtain a characterization of their face numbers. In particular we show that the convex hull of these face numbers is a cone. Finally we extend some of these results to arrangements of geodesics in surfaces of higher genus.
\end{abstract}
 \maketitle

\section*{Introduction}\label{intro}

Counting the number of connected components of a certain geometric set divided by its codimension-$1$ subsets is a classical problem in combinatorial geometry. The simplest possible (interesting) case is that of a partitioning of the Euclidean plane by finitely many straight lines. Such a collection determines a stratification of the plane consisting of \textit{vertices} (intersections of lines), \textit{edges} (maximal connected components of the lines not containing any vertex) and \textit{chambers} (maximal connected components of the plane containing neither the edges nor the vertices). The combinatorics that emerges from these intersections is intriguing. This is evident by the number of interesting problems and conjectures described in Gr\"unbaum's exposition \cite{gr72}. A systematic discussion of combinatorial aspects of hyperplane arrangements (i.e., higher-dimensional analogues of line arrangements) can be found in \cite[Chapter 18]{grunbaum67}.\par 

Classically, line arrangements are studied in the projective plane instead of the Euclidean plane. The first question that one can ask is to count the number of chambers formed by a line arrangement. An easy case (besides all concurrent lines) is that of lines in general position (i.e., no three lines are concurrent); here the number of chambers is $1 + \binom{n}{2}$, where $n$ is the number of lines. A formula for an arbitrary arrangement involves the M\"obius function of the intersection poset of the arrangement and it was discovered by Zaslavsky in \cite{zas75}. Note that Zaslavsky's theorem holds true in full generality for hyperplane arrangements. Let $f_i$ denote the number of $i$-dimensional strata, for $i = 0, 1, 2$, of the projective plane induced by a line arrangement. These are called as \textit{face numbers} and the triple $(f_0, f_1, f_2)$ is known as the \textit{$f$-vector} of a line arrangement. Many interesting questions arise when one wants to study relations between face numbers. For example, these numbers certainly satisfy the Euler relation $f_0 - f_1 + f_2 = 1$ but not all triples of natural numbers satisfying this relation arise as face numbers of line arrangements. One can find a list of known results and some conjectures in \cite[Section 2.2]{gr72}. In this paper we wish to answer similar questions but in the context of toric line arrangements. \par 

Partitioning problems for spaces other than Euclidean and projective spaces were studied by only handful of authors. To our knowledge the first paper that deals with a more general situation is by Zaslavsky \cite{zas77}. He derives a formula for counting the number of connected components a topological space when dissected by finitely many of its subspaces. He showed that not only the combinatorics of the intersections (which is encoded in the M\"obius function of the intersection poset) but also their geometry (as captured by the Euler characteristic) plays a role in determining the number of chambers. Pakula has considered arrangements of sub-spheres in a sphere in \cite{pakula93, pakula03}. In recent years several authors have considered toric arrangements. A toric arrangement is a finite collection of codimension-$1$ subtori in a torus. The formula for the number of chambers for such arrangements was first discovered by Ehrenborg et al. in \cite{ehr09}. The same formula was also independently discovered by Lawrence in \cite{Jim2011870} and by the second author in \cite{deshpande14}. Recently, Shnurnikov has characterized the set of all possible values of $f_2$ for toric line arrangements in \cite{shnurnikov11}. See also \cite{shnurnikov12} for arrangements in hyperbolic spaces, icosahedron and also arrangements of immersed circles in surfaces.\par 

The aim of this paper is to give a characterization of the $f$-vector for toric arrangements in a $2$-torus. The paper is organized as follows. In Section \ref{sec1} we introduce toric arrangements in full generality and fix notations. In Section \ref{sec2} we prove some properties of $f$-vectors for toric line arrangements. In particular we show that for toric line arrangements the convex hull of $f_0, f_2$ that appear as the number of vertices and as the number of chambers respectively is a cone in the first quadrant of the $(f_0, f_2)$-plane. Conversely, for every pair of integers in this cone there corresponds a toric line arrangement. Finally, in Section \ref{sec3} we outline future research by commenting on arrangements in surfaces of higher genus.

\section{Toric arrangements}\label{sec1}
The $l$-dimensional torus $\t^l$ is the quotient space $\R^l/\Z^l$. When identified with the set $[0, 1)^l$ it forms an abelian group with the group structure given by the componentwise addition modulo $1$. There is also a `multiplicative' way of looking at the torus when we consider it as the product of $S^1$'s. The group structure here is the componentwise multiplication of complex numbers of modulus $1$. However, throughout this paper, we stick to the additive way of looking at a torus. In this section we define toric arrangements and collect some relevant background material. \par 

We assume the reader's familiarity with basic algebraic topology and combinatorics. The combinatorics of posets and lattices that we need can be found in Stanley's book \cite[Chapter 3]{stan97}. As for the hyperplane arrangements Gr\"unbaum's book \cite[Chapter 18]{grunbaum67} covers mostly the enumerative aspect whereas the book of Orlik and Terao \cite{orlik92} describes modern results. The field of toric arrangements is fairly recent; Ehrenborg, Readdy and Slone mainly study the problem of enumerating faces of the induced decomposition of the torus in \cite{ehr09}. On the other hand a number theoretic aspect is explored by Lawrence in \cite{Jim2011870}. We also mention the pioneering work of De Concini and Procesi \cite{cpbook09}; they deal with aspects beyond the scope of this paper.\par 

We denote by $\pi\colon\thinspace\R^l\to\t^l$ the quotient map. Note that $\pi$ is also the covering map and $\R^l$ is the universal cover of the $l$-torus which is a compact manifold. We say that a $k$-subspace $V$ of $\R^l$ is \textit{rational} if it is the kernel of an $n\times l$ matrix $A$ with integer entries. The image $\ol{V} := \pi(V)$ is a closed subgroup of $\t^l$.
Topologically $\ol{V}$ is disconnected and each connected component is a $k$-torus. The connected components are known as \textit{toric subspaces} (or \textit{cosets}) of $\ol{V}$. Let $_0\ol{V}$ denote the coset containing $\mathbf{0}$ then $\ol{V}/_0\ol{V}$ is a finite abelian group whose order is the number of cosets of $\ol{V}$. One can check that every closed subgroup of the torus arises in this manner. It is important to note that the subgroup $\ol{V}$ depends only on the free abelian group generated by the row-space of $A$. Hence one can assume that the rows of $A$ form a basis for the row-space. The subgroup $\ol{V}$ is connected if and only if the greatest common divisor of all the $k\times k$ minors of $A$ is $1$. Two $k\times l$ matrices $A$ and $A'$ represent the same subgroup if and only if there exists a $k\times k$ unimodular matrix $U$ such that $A' = AU$.\par 

A \textit{toric hyperplane} is a toric subspace of codimension-$1$, i.e., it is the projection of an affine hyperplane in $\R^l$. We have the following definition. 

\bd{def1s1} A toric arrangement in $\t^l$ is a finite collection $\A = \{H_1,\dots, H_n\}$ of toric hyperplanes.\ed

A rational, codimension-$1$ subspace in $\R^l$ is specified by an equation $a_1 x_1 + \cdots + a_l x_l = c'$ where each $a_i\in\Z$. Hence we represent a toric hyperplane by a pair $(\mathbf{a}, c)$ where $\mathbf{a}$ is a row vector of integers and $c\in [0, 1)$. Consequently, sometimes it is convenient to express a toric arrangement as an augmented matrix $[A\mid \mathbf{c}]$ where $A$ is an $n\times l$ matrix of integers such that its each row represents the corresponding toric hyperplane and $\mathbf{c}$ is a vector in $[0, 1)^n$ representing intercept of each hyperplane.\par 

To every toric arrangement there is an associated periodic hyperplane arrangement $\tilde{\A}$ in $\R^l$. The inverse image of each $H_i$ under the covering map $\pi$ is the union of parallel integer translates of a codimension-$1$ subspace. Recall that a hyperplane arrangement is said to be \textit{essential} if the largest dimension of the subspace spanned by the normals to hyperplanes is $l$. We say that a toric arrangement is \textit{essential} if the associated hyperplane arrangement $\tilde{\A}$ is essential. Equivalently, it means that the rank of the matrix $A$ is $l$. Without loss of generality we assume that a toric arrangement is always essential; which forces $n\geq l$. If this is not the case then the enumerative problems that we consider in this paper reduce to equivalent problems in a torus of smaller dimension.\par 

The hyperplane arrangement $\tilde{\A}$ induces a stratification of $\R^l$ such that these open strata are relative interiors of convex polytopes. A nonempty  subset $F\subset \t^l$ is said to be a \textit{face} of the toric arrangement $\A$  if there is a strata $\tilde{F}$ of $\tilde{\A}$ such that $\pi(\tilde{F}) = F$. The \textit{dimension} of $F$ is the dimension of the support of $\tilde{F}$ and it is denoted by $\dim(F)$. Observe that every face of the toric arrangement lifts to a parallel class of strata in the periodic hyerplane arrangement $\tilde{\A}$. It is important to note that the closure of $F$ in $\t^l$ need not be homeomorphic to a disk. Hence a toric arrangement stratifies the ambient torus; this stratification need not define a regular cell structure but nonetheless has special properties. 

\bd{def2} A \textit{polytopal complex} is a cell complex $(X, \{e_{\lambda}\}_{\lambda\in\Lambda})$ with the following additional data.
\begin{enumerate}
\item Every cell $e_{\lambda}$ is equipped with a \textit{$k$-polytopal cell structure} which is a pair $(P_{\lambda}, \phi_{\lambda})$ of a $k$-convex polytope and a cellular map $\phi_{\lambda}\colon P_{\lambda}\to X$ such that $\phi_{\lambda}(P_{\lambda}) = \ol{e_{\lambda}}$ and the restriction of $\phi_{\lambda}$ to the interior of $P_{\lambda}$ is a homeomorphism.
\item If $e_{\mu}\cap \ol{e_{\lambda}}\neq \emptyset$ then $e_{\mu}\subset \ol{e_{\lambda}}$.
\item For every face $P'$ of $P_{\lambda}$, there exists a cell $e_{\mu}$ in $X$ and a map $b\colon Q_{\mu}\to \partial P_{\lambda}$ such that $b(\mathrm{Int} Q_{\mu}) = P'$ and $\phi_{\lambda}\circ b = \phi_{\mu}$.
\end{enumerate}\ed

The following lemma is a straightforward application of the fact that the covering map $\pi$ is stratification preserving. 

\bl{lem1} If $\A$ is a toric arrangement in $\t^l$ then the induced stratification is a polytopal complex. \el 

If the closure of each face is contractible then we say that the stratification defines a \textit{regular subdivision} of the torus or simply that it is a regular polytopal complex. A reason to elaborate on this type of cell structure is that we think it answers a question raised in \cite[Section 5]{ehr09}. One of the questions is about finding an analogue of regular subdivision of a manifold. The polytopal cell complex (or the totally normal cellularly stratified space as defined in \cite{tamaki01}) serves as the right analogue. The harder part of the question is the classification of flag $f$-vectors in this context. The combinatorial structure associated with this stratification is the following.

\bd{def3} Let $\A$ be a toric arrangement in $\t^l$. The \textit{face category} of $\A$, denoted by $\F(\A)$, is defined as follows. The objects of this category are faces of $\A$. A morphism from a face $\phi_{\mu}: P_{\mu}\to \ol{e_{\mu}}$ to another face $\phi_{\nu}: P_{\nu}\to \ol{e_{\nu}}$ is a map $b: P_{\mu}\to P_{\nu}$ such that $\phi_{\mu} = \phi_{\nu}\circ b$.
\ed

The face category is an acyclic category; which means that only the identity morphisms are invertible. The face category of a toric arrangement behaves much like the face poset of a hyperplane arrangement in the sense that the geometric realization of the category has the homotopy type of the torus (see \cite[Theorem 4.16]{tamaki01}). Moreover, if all the attaching maps are homeomorphisms (equivalently, closures of all faces are contractible) then the face category is equivalent to the underlying face poset \cite[Lemma 4.2]{tamaki01}. We refer the reader to \cite{tamaki01} for more on face categories. Now we move on to the next combinatorial object associated with a toric arrangement.

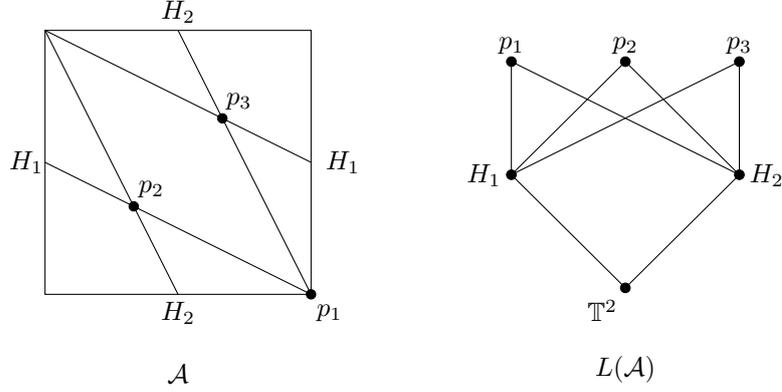
\begin{figure}[ht!]
\begin{center}
\begin{tikzpicture}[scale=3.5]
\draw[-] (0,0) rectangle ++(1,1);
\draw[-] (0,1) -- (0.5,0);
\draw[-] (0,1) -- (1,0.5);
\draw[-] (1,0) -- (0.5,1);
\draw[-] (1,0) -- (0,0.5);
\filldraw[-] (1/3,1/3) circle (0.5pt);
\filldraw[-] (2/3,2/3) circle (0.5pt);
\filldraw[-] (1,0) circle (0.5pt);
\node at (1.07,-0.07) {$p_1$};
\node at (0.4,0.4) {$p_2$};
\node at (0.73,0.73) {$p_3$};
\node at (0.5,-0.07) {$H_2$};
\node at (0.5,1.07) {$H_2$};
\node at (-0.07,0.5) {$H_1$};
\node at (1.12,0.5) {$H_1$};
\node at (0.5,-0.3) {$\A$};
\end{tikzpicture}
\quad\quad\quad
\begin{tikzpicture}[scale=3]
\filldraw[-] (1/2,0) circle (0.6pt);
\filldraw[-] (0,1/2) circle (0.6pt);
\filldraw[-] (0,1) circle (0.6pt);
\filldraw[-] (1/2,1) circle (0.6pt);
\filldraw[-] (1,1) circle (0.6pt);
\filldraw[-] (1,1/2) circle (0.6pt);
\draw[-] (1/2,0) -- (0,1/2);
\draw[-] (0,1/2) -- (0,1);
\draw[-] (0,1/2) -- (1/2,1);
\draw[-] (0,1/2) -- (1,1);
\draw[-] (1/2,0) -- (1,1/2);
\draw[-] (1,1/2) -- (1,1);
\draw[-] (1,1/2) -- (1/2,1);
\draw[-] (1,1/2) -- (0,1);
\node at (0.4,-0.1) {$\t^2$};
\node at (-0.12,1/2) {$H_1$};
\node at (1.12,1/2) {$H_2$};
\node at (0,1.07) {$p_1$};
\node at (1/2,1.07) {$p_2$};
\node at (1,1.07) {$p_3$};
\node at (0.5,-0.36) {$L(\A)$};
\end{tikzpicture}
\end{center}
\caption{A toric arrangement in $\t^2$.}
\label{fig1}
\end{figure}

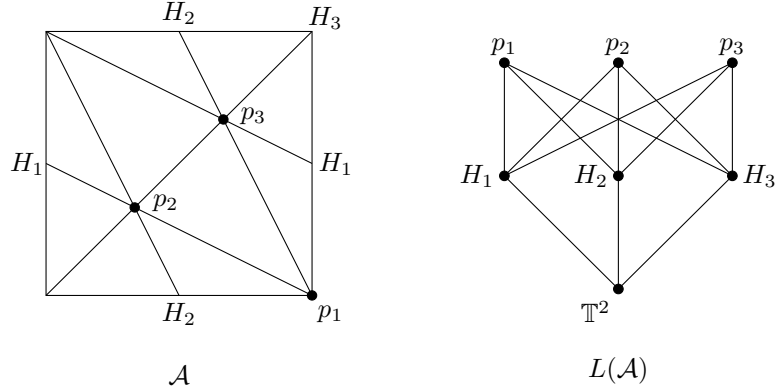
\begin{figure}[ht!]
\begin{center}
\begin{tikzpicture}[scale=3.5]
\draw[-] (0,0) rectangle ++(1,1);
\draw[-] (0,1) -- (0.5,0);
\draw[-] (0,1) -- (1,0.5);
\draw[-] (1,0) -- (0.5,1);
\draw[-] (1,0) -- (0,0.5);
\draw[-] (0,0) -- (1,1);
\filldraw[-] (1/3,1/3) circle (0.5pt);
\filldraw[-] (2/3,2/3) circle (0.5pt);
\filldraw[-] (1,0) circle (0.5pt);
\node at (1.07,-0.07) {$p_1$};
\node at (0.45,0.35) {$p_2$};
\node at (0.78,0.68) {$p_3$};
\node at (0.5,-0.07) {$H_2$};
\node at (0.5,1.07) {$H_2$};
\node at (-0.07,0.5) {$H_1$};
\node at (1.09,0.5) {$H_1$};
\node at (1.05,1.05) {$H_3$};
\node at (0.5,-0.3) {$\A$};
\end{tikzpicture}
\quad\quad\quad
\begin{tikzpicture}[scale=3]
\filldraw[-] (1/2,0) circle (0.6pt);
\filldraw[-] (0,1/2) circle (0.6pt);
\filldraw[-] (1/2,1/2) circle (0.6pt);
\filldraw[-] (0,1) circle (0.6pt);
\filldraw[-] (1/2,1) circle (0.6pt);
\filldraw[-] (1,1) circle (0.6pt);
\filldraw[-] (1,1/2) circle (0.6pt);
\draw[-] (1/2,0) -- (0,1/2);
\draw[-] (1/2,0) -- (1/2,1/2);
\draw[-] (1/2,0) -- (1,1/2);
\draw[-] (0,1/2) -- (0,1);
\draw[-] (0,1/2) -- (1/2,1);
\draw[-] (0,1/2) -- (1,1);
\draw[-] (1,1/2) -- (1,1);
\draw[-] (1,1/2) -- (1/2,1);
\draw[-] (1,1/2) -- (0,1);
\draw[-] (1/2,1/2) -- (1,1);
\draw[-] (1/2,1/2) -- (1/2,1);
\draw[-] (1/2,1/2) -- (0,1);
\node at (0.4,-0.1) {$\t^2$};
\node at (-0.12,1/2) {$H_1$};
\node at (0.38,1/2) {$H_2$};
\node at (1.12,1/2) {$H_3$};
\node at (0,1.07) {$p_1$};
\node at (1/2,1.07) {$p_2$};
\node at (1,1.07) {$p_3$};
\node at (0.5,-0.36) {$L(\A)$};
\end{tikzpicture}
\end{center}
\caption{A toric arrangement with regular cell decomposition.}
\label{fig2}
\end{figure}

\bd{def4}
The \textit{intersection poset} $L(\A)$ of a toric arrangement is defined to be the set of all connected components arising from all possible intersections of the toric hyperplanes ordered by reverse inclusion. By convention, the ambient torus corresponds to the empty intersection. The intersection poset is graded by the codimensions of the intersections. 
\ed
Before proceeding further let us look at a couple of examples. 

\be{ex1}
Let $\A$ be the toric arrangement in $\t^2$ obtained by projecting the lines $x = -2y$ and $y = -2x$. These toric hyperplanes intersect in three points $p_1 = (0, 0), p_2 = (1/3, 1/3)$ and $p_3 = (2/3, 2/3)$. The arrangement stratifies the torus into three $0$-faces, six $1$-faces and three $2$-faces. This is not a regular subdivision of the torus since the closure of every $2$-face is a cylinder. Figure \ref{fig1} shows the arrangement together with the associated intersection poset.
\ee

\be{ex2}
Now consider the arrangement formed by including the projection of the line $y = x$ in the previous arrangement. They intersect in the same three points as above. However, there are nine $1$-faces and six $2$-faces. The induced stratification is regular. Figure \ref{fig2} shows the arrangement and the associated intersection poset. 
\ee

Since our focus is on counting the number of various-dimensional faces of a toric arrangement we now turn to the combinatorics aspect. The idea that captures the combinatorics of the intersections is the M\"obius function of the arrangement which we now define. 

\bd{def5}
The M\"obius function of a toric arrangement is the function $\mu\colon L(\A)\times L(\A)\to \Z$ defined recursively as follows:
\[\ds\mu(X, Y) = \begin{cases}0, & \hbox{~if~} Y < X,\\ 1, &\hbox{~if~} X = Y, \\ -\sum_{X\leq Z< Y}\mu(X, Z), &\hbox{~if~} X < Y.\end{cases} \]
\ed
The M\"obius function plays an important role in counting the number of faces of an arrangement. The following theorem has appeared in \cite[Corollary 3.12]{ehr09}, \cite[Theorem 3]{Jim2011870} and \cite[Example 5.5]{deshpande14}.

\begin{theorem}\label{numkface}
Let $f_k$ denote the number of $k$-dimensional faces of a toric arrangement $\A$. Then we have 
\[ f_k  = \sum_{\substack{\dim Y = k\\ \dim Z = 0\\ Y\leq Z}} |\mu(Y, Z)|. \]
\end{theorem}

In particular the number of top-dimensional faces is determined by the values the M\"obius function takes at the points of intersections. The generating function for the face numbers is known as the \textit{$f$-polynomial} and defined as $f_{\A}(x) = \sum_{k=0}^l f_k x^{l-k}$.
Using Theorem \ref{numkface} above we get a particularly nice form for the $f$-polynomial
\[f_{\A}(x) = \sum_{\substack{\dim Y = k\\ \dim Z = 0\\ Y\leq Z}} |\mu(Y, Z)| x^{l - \dim Y}. \]
We say that the toric hyperplanes of an arrangement are in \textit{general position} if the intersection of any $i$ of the subtori, $i\geq 1$, is either empty or $(l-i)$-dimensional. A toric arrangement is called \textit{simple} if all the toric hypeprlanes are in general position. One can check that in case of simple arrangements every interval of the associated intersection poset is a Boolean algebra. For simple toric arrangements we have $f_{\A}(x) = f_0 (x + 1)^l$ hence 
\[f_k = f_0 \binom{l}{l-k}. \]
From the point of view of enumerative combinatorics simple arrangements are perhaps the easiest to understand.

\section{Some face enumeration formulas}\label{sec2}
In this section we focus our attention to arrangements in the $2$-torus $\t^2$ with the aim to explore relationship between the face numbers $f_0, f_1, f_2$ and $n$ the number of subtori in an arrangement. The projection of the straight line $ax+by=c$ in $\R^2$ under the canonical map $\pi$ onto the torus, where $a,b\in\Z$, is said to be a \textit{toric line}. Here we identify $c$ with $\pi(c)$ and assume that the $a$ and $b$ are coprime. Whenever convenient we will denote a line $l_j$ by an augmented matrix $[a_j, b_j\mid c_j]$ and say that the line is of type $(a_j, b_j)$ if the intercept is not relevant.\par 
The intersection two toric lines, of type say $(a_i, b_i)$ and $(a_j, b_j)$, is a finite set of points. The cardinality of the intersection is the absolute value of the determinant $\begin{vmatrix} a_i & b_i\\ a_j & b_j\end{vmatrix}$. The proof is a straightforward application of the Smith normal form (to be precise, structure theorem for finitely generated modules over PIDs). The Smith normal formal form of a $2\times 2$ matrix $A$ is a diagonal matrix with $1$ in the $(1, 1)$ position and $\det A$ in $(2, 2)$ position. The subgroup of the torus corresponding to $A$ is then a finite abelian group of order $|\det A|$ as the normal form is obtained by unimodular transformations (see \cite[Lemma 1]{shnurnikov11} for a geometric proof).

\bd{toriclar} A toric line arrangement is a finite collection $\A = \{l_1,\dots, l_n\}$ of toric lines in $\t^2$. \ed

As before we will denote $\A$ by an augmented matrix $[A\mid \mathbf{c}]$. We also assume that toric line arrangements are essential. Hence we do not consider the arrangements in which all lines are parallel.\par 

We now turn to the faces of an arrangement. For simplicity we call $0$-dimensional faces as vertices, $1$-dimensional faces as edges and $2$-dimensional faces as chambers; their numbers are denoted by $f_0, f_1, f_2$ respectively. These face numbers clearly satisfy the Euler relation $f_0 - f_1 + f_2 = 0$. It tells us that $f_1$ is redundant; hence we characterize pairs of natural numbers which appear as $(f_0, f_2)$ for some toric arrangement.\par

Let $L_0$ denote the set of all vertices. The number of lines in a toric line arrangement, that pass through a vertex $v$ is known as the \textit{degree} of that vertex and denoted by $\deg(v)$. The following is a straightforward application of Theorem \ref{numkface}.

\bl{degsum}\[f_1=\sum_{v\in L_0}\deg(v).\] \el


We denote by $t_j$ the number of vertices $v$ with $\deg v=j$. Since every vertex is formed by intersections, we must have $\deg v\ge2$ for each vertex $v$.

\bd{sgon}Let $\pi\colon\mathbb{R}^2\to T^2$ denote the canonical projection. A subset $C\subseteq T^2$ is a \textit{toric $k$-gon} if there exists an $k$-gon $C_0$ in $\R^2$ such that $\pi(C_0) = C$ and the restriction of $\pi$ to the interior of $C_0$ is a homeomorphism onto the image.
\ed

For a toric arrangement $\A$ let $p_k$ denote the number of chambers that are toric $k$-gons for $k\geq 3$. 

\bl{tj}The following results hold for any toric arrangement.
\begin{enumerate}
\item $\displaystyle f_0=\sum_jt_j$; \label{rel1}
\item $\displaystyle f_1=\sum_jjt_j$;\label{rel2}
\item $\displaystyle f_2=\sum_kp_k$;\label{rel3}
\item $\displaystyle f_2=\sum_j(j-1)t_j$.\label{rel4}
\end{enumerate} \el

\bpr Relations \ref{rel1} - \ref{rel3} follow from the definition and for \ref{rel4} one has to consider the Euler relation. \epr

\bl{kpk}Let $\A$ be a toric arrangement then we have that\[2f_1=\sum_kkp_k.\] \el

\bpr Given a chamber $C$ of $\A$ its lift $\pi^{-1}(C)$ consists of polygons in $\R^2$. The right hand side of the above equation is obtained by counting the edges bounding a polygon in that lift.  Since each such edge of this polygon projects downstairs to an edge in the arrangement, $\sum_{k\geq 3}kp_k$ counts every edge in the arrangement a certain number of times.\par 

Let an edge $e$ be counted $j$ times in the above manner. Since $e$ is arbitrary, if we show that $j=2$ we are done. Observe that a small enough neighbourhood $U$ of a point on $e$ is the union of $j$ semi-disks identified along their diameters. For a point $p\in e$ we have
\[H_2(\t^2,\t^2-p)\cong H_2(U,U-p)\cong H_1(U-p)\cong H_1(\bigvee_{j-1} S^1)\cong\Z^{j-1}.\]
Since $\t^2$ is a manifold, we already have $H_2(\t^2,\t^2-p)\cong\Z$ so that $j=2$.
\end{proof}

\bl{tjpkt2p3}
For any toric line arrangement, we have the following:
\begin{align}
\displaystyle t_2 &=\sum_{j\ge3}(j-3)t_j+\sum_{k\ge3}(k-3)p_k,\\
\displaystyle p_3 &=\sum_{j\ge2}2(j-2)t_j+\sum_{k\ge4}(k-4)p_k.
\end{align}
\el
\bpr The proof is a simple application of the Euler relation and definitions:

\item\[RHS-t_2=\sum_{j\ge2}(j-3)t_j+\sum_{k\ge3}(k-3)p_k=f_1-3f_0+2f_1-3f_2=0\] and
\item\[RHS-p_3=\sum_{j\ge2}2(j-2)t_j+\sum_{k\ge3}(k-4)p_k=2f_1-4f_0+2f_1-4f_2=0.\qedhere\]
\epr


We now turn our attention to Shnurnikov's result that characterizes the numbers that can occur as the number of chambers of a toric line arrangement. We reproduce the proof for the benefit of the reader.
  
\begin{theorem}[Shnurnikov {\cite[Theorem 1]{shnurnikov11}}]\label{shnut2f2}
Denote by $F(\t^2,n)$ the set of all possible values of $f_2$ that correspond to an arrangement of $n$ toric lines. Then \[F(\t^2,n)=\{n-1\}\cup\{l\in\N:l\ge 2n-4\}.\]
\et

\bpr
We first prove that $F(\t^2,n)\supseteq\{n-1\}\cup\{l\in\N:l\ge 2n-4\}$ by constructing arrangements with specified $f_2$. In order get $f_2=n-1$ consider the following arrangement:
\[\A = \{[1, 0\mid 0], [0, 1\mid c_j]: 1\leq j\leq n-1, c_i\neq c_j\hbox{~for~} i\neq j\}.\]

Figure \ref{fig3} below illustrates the construction for $n = 6$. The boundary edges of the fundamental domain correspond to toric lines $[1,0 \mid 0]$ and $[0,1\mid 0]$. Any of these lines are shown dotted if they are not part of the arrangement.
\begin{figure}[ht]
\hfill
\begin{tikzpicture}[scale=3]
\draw[dashed] (0,0) -- (1,0);
\draw[dashed] (0,1) -- (1,1);
\draw[dashed] (0,0) -- (0,1);
\draw[dashed] (1,0) -- (1,1);
\draw[-] (0.6,0) -- (0.6,1);
\draw[-] (0.2,0) -- (0.2,1);
\draw[-] (0.3,0) -- (0.3,1);
\draw[-] (0.4,0) -- (0.4,1);
\draw[-] (0.5,0) -- (0.5,1);
\draw[-] (0,0.5) -- (1,0.5);
\end{tikzpicture}
\hfill\text{}\caption{$(n,f_2)=(6,5)$}\label{fig3}
\end{figure}
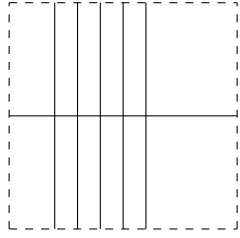

We now construct an arrangement $\A$ with $f_2=2n-4+a$ for an arbitrary whole number $a$. Consider the arrangement consisting of the toric lines of the following types:
\begin{enumerate}
\item $[0,1\mid 0]$;
\item $[a+1,-1\mid 0]$; 
\item $[1,0\mid 0]$, $\ds\left[1,0\mid\frac1{a+2}\right]$,$\ds\left[1,0\mid\frac1{a+3}\right],\dots, \ds\left[1,0\mid\frac1{a+n-2}\right]$.
\end{enumerate}
The construction is illustrated in Figure \ref{fig4} for $n = 6$.\par 
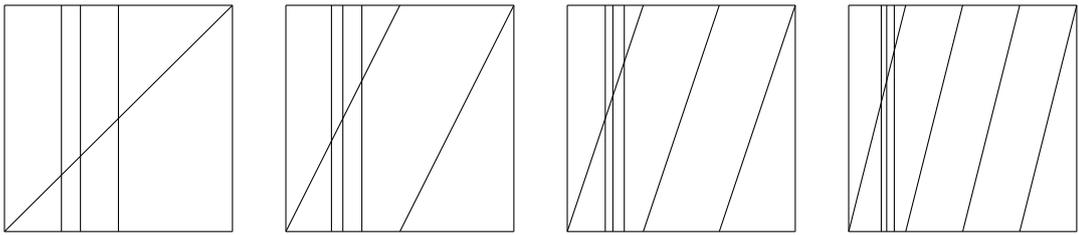
\begin{figure}[ht]
\hfill\begin{tikzpicture}[scale=3]
\draw[-] (0,0) -- (1,0);
\draw[-] (0,1) -- (1,1);
\draw[-] (0,0) -- (0,1);
\draw[-] (1,0) -- (1,1);
\draw[-] (0,0) -- (1,1);
\draw[-] (1/2,0) -- (1/2,1);
\draw[-] (1/3,0) -- (1/3,1);
\draw[-] (1/4,0) -- (1/4,1);
\end{tikzpicture}\hfill
\begin{tikzpicture}[scale=3]
\draw[-] (0,0) -- (1,0);
\draw[-] (0,1) -- (1,1);
\draw[-] (0,0) -- (0,1);
\draw[-] (1,0) -- (1,1);
\draw[-] (0,0) -- (1/2,1);
\draw[-] (1/2,0) -- (1,1);
\draw[-] (1/3,0) -- (1/3,1);
\draw[-] (1/4,0) -- (1/4,1);
\draw[-] (1/5,0) -- (1/5,1);
\end{tikzpicture}\hfill
\begin{tikzpicture}[scale=3]
\draw[-] (0,0) -- (1,0);
\draw[-] (0,1) -- (1,1);
\draw[-] (0,0) -- (0,1);
\draw[-] (1,0) -- (1,1);
\draw[-] (0,0) -- (1/3,1);
\draw[-] (1/3,0) -- (2/3,1);
\draw[-] (2/3,0) -- (1,1);
\draw[-] (1/4,0) -- (1/4,1);
\draw[-] (1/5,0) -- (1/5,1);
\draw[-] (1/6,0) -- (1/6,1);
\end{tikzpicture}\hfill
\begin{tikzpicture}[scale=3]
\draw[-] (0,0) -- (1,0);
\draw[-] (0,1) -- (1,1);
\draw[-] (0,0) -- (0,1);
\draw[-] (1,0) -- (1,1);
\draw[-] (0,0) -- (1/4,1);
\draw[-] (1/4,0) -- (1/2,1);
\draw[-] (1/2,0) -- (3/4,1);
\draw[-] (3/4,0) -- (1,1);
\draw[-] (1/5,0) -- (1/5,1);
\draw[-] (1/6,0) -- (1/6,1);
\draw[-] (1/7,0) -- (1/7,1);
\end{tikzpicture}\hfill\text{}\caption{$(n,f_2)=(6,8),(6,9),(6,10)$ and $(6,11)$}\label{fig4}
\end{figure}
In order to prove the reverse containment assume that $\A$ contains at most $m$ toric lines of the same type (i.e., $m$ is the maximal number of parallel lines). If $m = n-1$ then $f_2$ is a multiple of $n-1$. To see this observe that the $n$th line intersects all the previous lines in the same number of points. \par
For $2\leq m\leq n-2$, each of the remaining $n-m$ lines intersect the $m$ parallel lines in at least $m$ points. Hence we have $f_2 \geq m (n-m)$. An easy exercise in calculus shows that the function $x(n-x)$ attains its bounds on the interval $[2, n-2]$ and the minima is $2(n-2)$ which implies that 
\[f_2 \geq 2n - 4. \]
The last case is that of having no parallel lines in $\A$. Here one has to consider several sub-cases. First, assume that any two of the lines intersect in at least two points. So, if $n = 2$ then $f_2\geq 2$. By induction on $n$ assume that for $n-1$ lines $f_2 \geq 2(n -1) - 2$. If the $n$th line creates a new vertex then $t_2$ goes up by $1$ and if it passes through an existing vertex then $t_j$ goes down by one and $t_{j+1}$ goes up by one for some $j$. In either case, using Identity \ref{rel4} of Lemma \ref{tj}, we get that $f_2$ increases by at least 2.\par 
Assume that some two lines, say $l_i, l_j$, meet in exactly one point. If all the remaining $n-2$ lines meet $l_i\cup l_j$ in at least two points then $f_2 \geq 2n - 3$. Otherwise we claim that there exists a third line $l_k$ which passes through the same intersection point. In order to prove the claim assume that $l_i$ is of type $(1, 0)$ and $l_j$ is of type $(0, 1)$ then as the line $l_k$ should intersect both these lines in point it has to be of the type $(1, 1)$ or $(1, -1)$. In either case the intersection $l_i\cap l_j\cap l_k$ is singleton. Furthermore it is easy to prove that if $l'$ is any line which is not parallel to either $l_i, l_j$ or $l_k$ then it intersects these three lines in at least two points.  In this case $f_2 \geq 2n -4$; the proof is on the same lines as that of the first sub-case.
\epr


It is well-known that the face numbers $f_0, f_2$ of projective line arrangements satisfy linear inequalities. These inequalities are such that their convex hull is a cone in $(f_0, f_2)$-plane. However not all lattice point in that cone are realizable as face numbers of projective line arrangements (see \cite[page 401]{grunbaum67} for details). On the other hand face numbers of convex polyhedra also satisfy similar inequalities and also determine a cone. Interestingly every pair $(f_0, f_2)$ satisfying these inequalities indeed corresponds to some convex polyhedron (see \cite[page 190]{grunbaum67} for details). The case of toric arrangements is not very different as proved below. We say that a toric line arrangement is \textit{simplicial} if all the chambers are triangles. 

\bt{f0f2} Given $f_0,f_2\in\N$, there exists a toric arrangement $\A$ with $f_0$ vertices and $f_2$ faces if and only if \[f_0\le f_2\le2f_0.\] Equality on the left holds if and only if $\A$ is simple; equality on the right holds if and only if $\A$ is simplicial.\et

\bpr
We prove the `only if' part first. Using Lemma \ref{tj} we see that:
\[f_0=\sum_{j\ge2}t_j\le\sum_{j\ge2}(j-1)t_j=f_2.\]
The second inequality can be written as $\ds f_2\le 2f_1-2f_2$, or equivalently as $\ds 2f_1\ge3f_2$ which follows from Lemmas \ref{tj} and \ref{kpk}:
\[2f_1=\sum_{k\ge3}kp_k\ge\sum_{k\ge3}3p_k=3f_2.\]

The `if' part on the other hand can be proved constructively. Consider the arrangement $\A$ which contains the following $n = f_2 - f_0 + 2$ toric lines : 
\begin{enumerate}
\item $[0,1\mid 0]$;
\item $[f_0,-1\mid 0]$;
\item $\ds\left[1,0\mid \frac r{f_0}\right]$ for all $0\le r<f_2-f_0$.
\end{enumerate}
See Figures \ref{fig5}, \ref{fig6} for illustrations. \par 
\begin{figure}[ht!]
\begin{tikzpicture}[scale=0.8]
\draw[-] (0,0) node[below] {$v_1$} -- (1,4) node[above] {$v_2$};
\draw[-] (1,0) node[below] {$v_2$} -- (2,4) node[above] {$v_3$};
\draw[-] (2,0) node[below] {$v_3$} -- (3,4) node[above] {$v_4$};
\draw[-] (3,0) node[below] {$v_4$} -- (4,4);
\draw[-] (0,0) -- (4,0);
\draw[-] (0,4) -- (4,4);
\draw[dashed] (0,0) -- (0,4) node[above] {$v_1$};
\draw[dashed] (4,0) node[below] {$v_1$} -- (4,4) node[above] {$v_1$};
\end{tikzpicture}
\begin{tikzpicture}[scale=0.8]
\draw[-] (0,0) node[below] {$v_1$} -- (1,4) node[above] {$v_2$};
\draw[-] (1,0) node[below] {$v_2$} -- (2,4) node[above] {$v_3$};
\draw[-] (2,0) node[below] {$v_3$} -- (3,4) node[above] {$v_4$};
\draw[-] (3,0) node[below] {$v_4$} -- (4,4);
\draw[-] (0,0) -- (4,0);
\draw[-] (0,4) -- (4,4);
\draw[-] (0,0) -- (0,4) node[above] {$v_1$};
\draw[-] (4,0) node[below] {$v_1$} -- (4,4) node[above] {$v_1$};
\end{tikzpicture}
\begin{tikzpicture}[scale=0.8]
\draw[-] (0,0) node[below] {$v_1$} -- (1,4) node[above] {$v_2$};
\draw[-] (1,0) node[below] {$v_2$} -- (2,4) node[above] {$v_3$};
\draw[-] (2,0) node[below] {$v_3$} -- (3,4) node[above] {$v_4$};
\draw[-] (3,0) node[below] {$v_4$} -- (4,4);
\draw[-] (0,0) -- (4,0);
\draw[-] (0,4) -- (4,4);
\draw[-] (0,0) -- (0,4) node[above] {$v_1$};
\draw[-] (4,0) node[below] {$v_1$} -- (4,4) node[above] {$v_1$};
\draw[-] (1,0) -- (1,4);
\end{tikzpicture}\caption{From left: toric arrangements with $f$-vectors (4,8,4), (4,9,5) and (4,10,6)}
\label{fig5}
\end{figure}
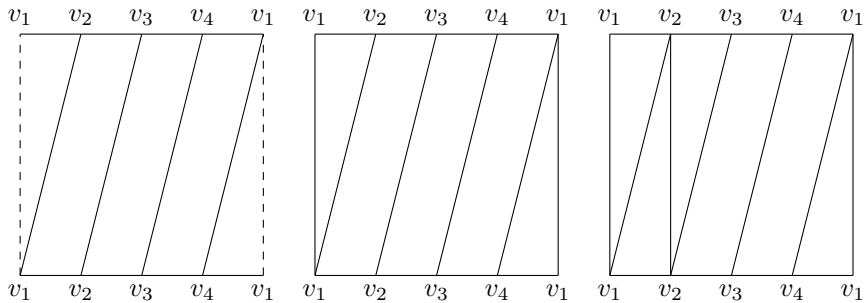

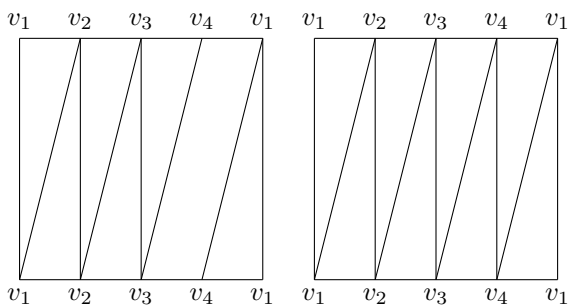
\begin{figure}[ht!]
\begin{tikzpicture}[scale=0.8]
\draw[-] (0,0) node[below] {$v_1$} -- (1,4) node[above] {$v_2$};
\draw[-] (1,0) node[below] {$v_2$} -- (2,4) node[above] {$v_3$};
\draw[-] (2,0) node[below] {$v_3$} -- (3,4) node[above] {$v_4$};
\draw[-] (3,0) node[below] {$v_4$} -- (4,4);
\draw[-] (0,0) -- (4,0);
\draw[-] (0,4) -- (4,4);
\draw[-] (0,0) -- (0,4) node[above] {$v_1$};
\draw[-] (4,0) node[below] {$v_1$} -- (4,4) node[above] {$v_1$};
\draw[-] (1,0) -- (1,4);
\draw[-] (2,0) -- (2,4);
\end{tikzpicture}
\begin{tikzpicture}[scale=0.8]
\draw[-] (0,0) node[below] {$v_1$} -- (1,4) node[above] {$v_2$};
\draw[-] (1,0) node[below] {$v_2$} -- (2,4) node[above] {$v_3$};
\draw[-] (2,0) node[below] {$v_3$} -- (3,4) node[above] {$v_4$};
\draw[-] (3,0) node[below] {$v_4$} -- (4,4);
\draw[-] (0,0) -- (4,0);
\draw[-] (0,4) -- (4,4);
\draw[-] (0,0) -- (0,4) node[above] {$v_1$};
\draw[-] (4,0) node[below] {$v_1$} -- (4,4) node[above] {$v_1$};
\draw[-] (1,0) -- (1,4);
\draw[-] (2,0) -- (2,4);
\draw[-] (3,0) -- (3,4);
\end{tikzpicture}\caption{From left: toric arrangements with $f$-vectors (4,11,7) and (4,12,8)}\label{fig6}
\end{figure}
Now assume that $\A$ is an arrangement with $f_2 = f_0$. Then $\sum_{j\geq 2}(j-1)t_j = \sum_{j\geq 2} t_j$ implies that 
\[\sum_{j\geq 2} (j -2) t_j = 0. \]
Since $(j-2) > 0$ for every $j \geq 3$ we have that $t_j = 0$ for those $j$'s. Consequently there are only degree $2$ vertices; equivalently the arrangement is simple. Converse of this statement is also clear.\par 
Now assume that $f_2 = 2f_0$. Using the equation 
\[\sum_{k\geq 3}p_k = 2 \sum_{k\geq 3}(\frac{k}{2} - 1)p_k\]
we get that $p_k = 0$ for $k\geq 4$. Converse can be proved analogously. \epr

Combining above inequality with Theorem \ref{shnut2f2} we see that for an arrangement of $n$ toric lines $f_0 \in \{\lfloor\frac{n-1}{2}\rfloor \}\cup\{l : l \geq  n- 2\}$. As $f_0$ is not bounded above there is no hope for a complete characterization of the pairs $(n, f_0)$. Instead we focus on the triples $(n,f_0,f_2)$. More precisely we would like to characterize all such triples of natural numbers for which there exists an arrangement $\A$ of $n$ lines, with $f_0$ vertices and $f_2$ faces. For $n \geq 2$ let 
\[ \mathcal{C}'(n) := \{(f_0, f_2) \mid f_0\leq f_2\leq 2f_0, f_2\geq n-1\}\]
we call it the \textit{potential search region} for toric arrangements of $n$ lines. There is an obvious chain of inclusions $\mathcal{C}'(2) \supset \cdots \mathcal{C}'(n) \supset \mathcal{C}'(n+1)\supset \cdots$. Our aim is to characterize elements of $\mathcal{C}'(n)$ that are realizable as face numbers of toric arrangement we denote this subset by $\mathcal{C}(n)$. We start with the easiest case. 

\bl{2f0f2}
For toric arrangements of $2$ lines we have 
\[\mathcal{C}(2) = \{(f, f)\mid f\in \mathbb{N} \}. \]
Equivalently, all toric arrangements of $2$ lines are simple and their $f$-vectors are of the form $(f, 2f, f)$ for all natural numbers $f$. \el
\bpr Since there are only 2 lines, all vertices have degree 2. This ensures, $f_0=t_2=f_2$. On the other hand, if $f_0=f_2$ is given, then indeed consider $\A$ with two toric lines, one each of the types $(0 , 1)$ and $(f_0 , -1)$.\epr

The next case, i.e. complete description of $\mathcal{C}(3)$ is difficult. First, just like in the $2$ lines case, the points of the type $(f_0, f_0)$ are completely realizable using simple arrangements. Second, we show that not all integer points of the type $(f_0, 2f_0)$ are realizable. 
\bt{3f0eqf2} There exists an arrangement of 3 lines with $f_0=f_2$ if and only if $f_0\ge2$ \et

\bpr The `only if' part follows at once from Theorem \ref{shnut2f2}.\par 
For the `if' part set $\A=\{[1,0\mid 0],[0,1\mid 0],[f_0-2,-1\mid r]\}$ where $r$ is a fixed irrational number. The irrationality of $r$ ensures that the line $[f_0-2,-1\mid r]$ does not pass through the intersection of the other two lines. The reader can easily check that there are $f_0$ vertices, all of which have degree $2$, so that $f_2=f_1-f_0=2t_2-t_2=t_2=f_0$. \epr

\bt{3f0f2} There exists an arrangement of 3 lines with $f_0$ vertices and $f_2=2f_0$ chambers if and only if $f_0$ is an odd number \et
\bpr Start by assuming that $f_0$ is even. Since there are only 3 lines, it is clear that 
\[f_2 = t_2 + 2t_3.\]
As $2f_0 = 2t_2 + 2t_3$ we have $t_2=0$ and consequently all vertices are of degree 3. Whence all the three lines $l_1, l_2, l_3$ in the arrangement pass through all the vertices. This shows that any two of the $l_i$'s intersect at $f_0$ many points. Without loss of generality assume that $l_1$ is of type $(0, 1)$ and $l_2$ of type $(f_0, -a)$, where of course $a, f_0$ are coprime. Since $l_3$ intersects $l_1$ as well as $l_2$ at $f_0$ vertices, $l_3$ must be of the type $(f_0, -a\pm 1)$.\par 

Now if $f_0$ is even then $a$ must be an odd  number coprime to $f_0$. This means that $-a\pm1$ is even and hence the line of type $(f_0 , -a\pm1)$ has two components, which is a contradiction, for it gives an arrangement of 4 lines.\par 

Conversely, if $f_0$ is odd, say $2k-1$, consider the arrangement of 3 lines one each of the types $(k,-(k-1))$, $(k-1,-k)$ and $(1,1)$. This gives an arrangement with $2k-1$ vertices and $4k-2$ faces.\epr

Now we prove a necessary condition for the points of the type $f_0 < f_2 < 2f_0$ to be realizable. 
\bt{thmdivf0} For an arrangement of $3$ toric lines such that $f_0 < f_2 < 2f_0$ then $f_2 - f_0$ divides $f_0$. \et 
\bpr The proof is straightforward once the reader realizes that it is enough to show that $t_3\mid f_0$. Hence we leave it as a simple exercise for the reader. \epr

The above theorem implies, for example, that there can not be an arrangement of three toric lines such that $f_0 = 4$ and $f_2 = 7$. We do not claim that this is also a sufficient condition. In general, complete characterization of $\mathcal{C}(n)$ seems to be a hard problem. We end this section by a result that characterizes arrangements for which the degree of the vertices is constant.

\bp{lem1sec2} Let $\A$ be a toric arrangement of $n$ lines such that $\deg(v) = k, \forall v\in L_0$ then $\A$ is either simple or simplicial.\ep
\bpr We have that $t_k\neq 0$ for some $k\geq 2$ and all other $t_i$'s are zero. Therefore,
\[f_2 = (k-1) t_k = (k-1) f_0 \leq 2 f_0. \]
Thus $k$ is either $2$ or $3$. The $k = 2$ case implies that no three lines are concurrent which means that the arrangement is simple. Whereas as for $k = 3$ using Lemma \ref{tjpkt2p3} we see that $p_k = 0$ for $k \geq 4$.\epr

\section{Concluding Remarks} \label{sec3}
We end the paper by a brief discussion about possible directions for future research. One direction is to look at arrangements in surfaces of higher genus and the other direction is to study these problems in higher-dimensional tori. In \cite[\S 6]{shnurnikov12}, Shnurnikov has realized the genus $g$ surface $M_g$ as the quotient of the hyperbolic plane $\h^2$ by a certain discrete subgroup $G$ of $\mathrm{Isom}(\h^2)$. A \textit{simple closed geodesic} in $M_g$ is defined to be the image of a geodesic line in $\h^2$ under the covering projection $\h^2\rightarrow\h^2/G$

\bd{arrgenus}
A finite collection of simple closed geodesics in $M_g$ inducing a polytopal cell structure is known as a \textit{geodesic arrangement} in a genus-$g$ surface.
\ed

We analogously define the intersection poset $L(\A)$ and the face numbers $f_0, f_1, f_2$ for geodesic arrangements. Some of the results proved in Section \ref{sec2} easily generalize in this case. They remain unchanged except for an additional term $\chi(M_g)$ the Euler characteristic. We enumerate such results now.

\bt{numkfacemg} For a geodesic arrangement $\A$ in $M_g$ we have\[f_1=\sum_{v\in L_0}\deg v.\] \et

\bpr By \cite[Theorem 4.2]{deshpande14}:
\[f_1=\sum_{\dim Y=1}\left(-\sum_{\substack{Z\in L(\A)\\Y\le Z}}\mu(Y,Z)\chi(Z)\right).\]
Since $\chi(Y)=0$ we could replace $Y\le Z$ with $Y<Z$ above. The summation then is over all $Z$ that cover $Y$, so that all $\mu$ values are -1. So the summation becomes:
\[f_1=\sum_{\dim Y=1}\sum_{\substack{Z\in L(\A)\\Y<Z}}\chi(Z).\]
Finally, $\chi(Z)=1$ for the $Z$ are just points.
\[f_1=\sum_{\substack{\dim Z=0\\Z\in L(\A)}}\left(\sum_{\substack{\dim Y=1\\Y<Z}}1\right).\]
But the inner summation is the degree of the vertex $Z$ and we are done.
\end{proof}

Since we are proving analogous results about vertices, degrees as well as $k$-gons, we now clarify the notion of an $s$-gon in a geodesic arrangement.

\bd{sgonmg}Let $\pi\colon\h^2\to M_g$ denote the covering map. A subset $C\subseteq M_g$ is a $k$-gon if there exists a geodesic $k$-gon $C_0$ of $\h^2$ such that $\pi(C_0) = C$ and restriction of $\pi$ to the interior of $C_0$ is a homeomorphism onto the image. \ed
As before $t_j$ stands for the number of degree $j$ vertices and $p_k$ stands for the number of $k$-gons.
\bl{jtjmg}
The following results hold for a geodesic arrangement on a genus $g$ surface:
\begin{enumerate}
\item $\ds f_0=\sum_jt_j$;
\item $\ds f_1=\sum_jjt_j$;
\item $\ds f_2=\sum_kp_k$;
\item $\ds2f_1=\sum_kkp_k$;
\item $\ds f_2-\chi(M_g)=\sum_j(j-1)t_j$;
\item $\ds2(f_0-\chi(M_g))=\sum_k(k-2)p_k$.
\end{enumerate}
\el
\bpr
(1),(2) and (3) are by definition. (4) is proved analogously as in \S\ref{sec2}.
(5) follows from (1),(2) and the Euler relation. (6) is a consequence of (3),(4) and the Euler relation. \end{proof}

\bl{tjpkt2p3mg}
For a geodesic arrangement in a genus-$g$ surface, we have the following:
\begin{align*}
\displaystyle t_2-3\chi(M_g)=\sum_{j\ge3}(j-3)t_j+\sum_{k\ge3}(k-3)p_k, \\
\displaystyle p_3-4\chi(M_g)=\sum_{j\ge2}2(j-2)t_j+\sum_{k\ge4}(k-4)p_k.
\end{align*} \el
\bpr The proof follows from the application of the Euler relation and Lemma \ref{jtjmg} above. \end{proof}
We now state a partial analogue of Theorem \ref{f0f2}; the proof is on the similar lines.
\bt{f0f2mg}
For an arrangement on a genus-$g$ surface, we have the following:\[f_0+\chi(M_g)\le f_2\le2(f_0-\chi(M_g)).\]
The equality on the left holds if and only if the arrangement is simple (i.e., $t_j = 0$ for $j\geq 3$) whereas the equality on the right holds if and only if the arrangement is simplicial (i.e., $p_k = 0$ for $k\geq 4$). \et

However the complete characterisation of tuples $(f_0,f_2)$ - for which there is a geodesic arrangement in a genus-$g$ surface with $f_0$ vertices and $f_2$ faces - analogous to Theorem \ref{f0f2} is a work in progress.\par 

As for the toric arrangements in higher-dimensional tori we refer the reader to the recent work of Shnurnikov \cite{shnurnikovarxiv} where the numbers that appear as number of chambers of toric arrangement are characterized.

\bibliographystyle{abbrv} 
\bibliography{toricref} 

\end{document}